\documentclass[12pt]{amsart}
\usepackage{amssymb}
\usepackage{amsmath,amscd}

\let\frak\mathfrak
\let\Bbb\mathbb

\def\>{\relax\ifmmode\mskip.666667\thinmuskip\relax\else\kern.111111em\fi}
\def\<{\relax\ifmmode\mskip-.333333\thinmuskip\relax\else\kern-.0555556em\fi}
\def\vsk#1>{\vskip#1\baselineskip}
\def\vv#1>{\vadjust{\vsk#1>}\ignorespaces}
\def\vvn#1>{\vadjust{\nobreak\vsk#1>\nobreak}\ignorespaces}

\let\Medskip\medskip
\def\medskip{\par\Medskip}
\let\Bigskip\bigskip
\def\bigskip{\par\Bigskip}

\let\Maketitle\maketitle
\def\maketitle{\hrule height0pt\vskip-\baselineskip
\Maketitle\thispagestyle{empty}\let\maketitle\empty}

\numberwithin{equation}{section}

\theoremstyle{definition}

\let\mc\mathcal
\let\nc\newcommand

\nc{\on}{\operatorname}
\nc{\Z}{{\mathbb Z}}
\nc{\C}{{\mathbb C}}
\nc{\N}{{\mathbb N}}
\nc{\pone}{{\mathbb C}{\mathbb P}^1}
\nc{\arr}{\rightarrow}
\nc{\larr}{\longrightarrow}
\nc{\al}{\alpha}
\nc{\W}{{\mc W}}
\nc{\la}{\lambda}
\nc{\su}{\widehat{{\mathfrak sl}}_2}
\nc{\g}{{\mathfrak g}}
\nc{\h}{{\mathfrak h}}
\nc{\m}{{\mathfrak m}}
\nc{\n}{{\mathfrak n}}
\nc{\Gm}{\Gamma}
\nc{\La}{\Lambda}
\nc{\gl}{\widehat{\mathfrak{gl}_2}}
\nc{\bi}{\bibitem}
\nc{\om}{\omega}
\nc{\Res}{\on{Res}}
\nc{\gm}{\gamma}
\nc{\Om}{\Omega}

\def\Res{\on{Res}}

\def\B{{\mc B}}

\def\V{{\mc V}}

\let\Dl\Delta

\let\ge\geqslant
\let\geq\geqslant

\nc{\gln}{\mathfrak{gl}_N}
\nc{\sln}{\mathfrak{sl}_N}

\def\beq{\begin{equation}}
\def\eeq{\end{equation}}
\def\be{\begin{equation*}}
\def\ee{\end{equation*}}

\nc{\bean}{\begin{eqnarray}}
\nc{\eean}{\end{eqnarray}}
\nc{\bea}{\begin{eqnarray*}}
\nc{\eea}{\end{eqnarray*}}
\nc{\bs}{\boldsymbol}
\nc{\Ref}[1]{{\rm(\ref{#1})}}
\nc{\glN}{\mathfrak{gl}_N}
\nc{\glNt}{\mathfrak{gl}_N[t]}
\nc{\s}{sing}
\nc{\R}{\Bbb R}

\nc{\Oml}{{\Om_{\bs\la}}}
\nc{\OmLb}{{\Om_{\bs\La,\bs\la,\bs b}}}
\nc{\Ol}{{\mc O_{\bs\la}}}
\nc{\OLb}{{\mc O_{\bs\La,\bs\la,\bs b}}}
\nc{\VSl}{{(\V^S)_{\bs\la}}}
\nc{\Bl}{{\B_{\bs\la}}}
\nc{\Ml}{{\mc M_{\bs\la}}}
\nc{\Mlb}{{\mc M_{\bs\La,\bs\la,\bs b}}}
\nc{\Blb}{{\B_{\bs\La,\bs\la,\bs b}}}
\nc{\Omn}{{\Omega_{\bs n,\bs b,\bs K}}}
\nc{\Omlb}{{\bar\Om_{\bs\la}}}
\nc{\ep}{\epsilon}

\nc{\Dlb}{\Dl_{\bs\La,\bs\la,\bs b,\bs K}}
\nc{\Bb}{{\bf b}}
\nc{\glt}{{\frak{gl}_2}}
\nc{\A}{{\mc A}}
\nc{\slt}{{\frak{sl}_2}}
\nc{\Ma}{{\mc M_{\bs a}}}
\nc{\Mal}{{\mc M_{\bs\la,\bs a}}}
\nc{\Malp}{{\mc M_{\phi,\bs\la,\bs a}}}

\nc{\Bal}{{\B_{\bs\la,\bs a}}}
\nc{\Ola}{{\mc O_{\bs\la,\bs a}}}
\nc{\Bv}{{\mc B_{\V^S}}}
\nc{\Bvz}{{\mc B^0_{\V^S}}}

\nc{\sing}{{\rm Sing\,}}
\nc{\Uglt}{U(\glt)}
\nc{\Olo}{{\mc O^0_{\bs\la}}}
\nc{\kk}{K}

\nc{\Oll}{{\Omega_{\bs\la}}}
\nc{\T}{{\mc T}}

\nc{\CC}{{\mc C}}
\nc\Vl{{(\V^S)^{sing}_{\bs\la}}}

\nc{\PP}{{\Bbb P}}
\nc{\LL}{{\mc L}}
\nc{\FF}{{\mc F}}

\begin{document}

\nc{\codim}{\on{codim}}
\nc{\supp}{\on{supp}}

\title[Intersection cohomology of the complement of an arrangement]
{Intersection cohomology of a rank one local system on the complement of a hyperplane-like divisor}

\author[D.\,Arinkin, A.\,Varchenko]
{D.\,Arinkin$^\dagger$,  A.\,Varchenko$^\diamond$}

\maketitle

\begin{center}
{\it Department of Mathematics, University of North Carolina
at Chapel Hill\\ Chapel Hill, NC 27599-3250, USA\/}
\end{center}

{\let\thefootnote\relax
\footnotetext{\vsk-.8>\noindent
$^\dagger$\,Supported in part by the Sloan Fellowship}}
{\let\thefootnote\relax
\footnotetext{\vsk-.8>\noindent
$^\diamond$\,Supported in part by NSF grant DMS-1101508}}

\medskip

\begin{abstract}
Under a certain condition A we give a construction to calculate the intersection cohomology of
a rank one local system on the complement to a hyperplane-like divisor.

\end{abstract}

\maketitle

\bigskip

Let $X$ be a smooth connected complex manifold and
$D$ a divisor. The divisor $D$ is {\it hyperplane-like}
if $X$ can be covered by coordinate charts such that in each chart
$D$ is the union of
hyperplanes. Such charts  are called \emph{linearizing}.

Let $D$ be a hyperplane-like divisor, $V$
a linearizing chart. A \emph{local edge} of $D$ in $V$ is
any nonempty irreducible intersection in $V$ of hyperplanes of $D$ in $V$.
A local edge  is \emph{dense}  if the subarrangement of all hyperplanes in $V$ containing
the edge is irreducible: the hyperplanes cannot be partitioned into nonempty
sets so that, after a change of coordinates, hyperplanes in different
sets are in different coordinates.
An \emph{edge} of $D$ is the maximal analytic continuation in $X$ of a local edge.
Any edge is an immersed submanifold
in $X$. An edge is called dense if it is locally dense.

\medskip

Let $U=X-D$ be the complement to $D$. Let $\LL$ be a rank one local system
on $U$ with nontrivial monodromy around each irreducible component of $D$.
We want to calculate the intersection cohomology $H^*(X; j_{!*}\LL)$ where $j:U\to X$ is the embedding.
In this calculation a role will be played by the dual local system $\LL^\vee$ on $U$ with the inverse monodromy.

Let $\pi : \tilde X\to X$ be a resolution of singularities of $D$ in $X$.

\medskip
\noindent
{\bf Remark.}\ Note that among the resolutions of singularities
there is the minimal one.
The minimal resolution is constructed by first blowing-up dense vertices of $D$,
then by blowing-up the proper preimages of dense
one-dimensional edges of $D$ and so on, see \cite{V, STV}.

The preimage $\pi^{-1}(D)$ is a hyperplane-like divisor in $\tilde X$.
 The inverse image $\pi^*\LL$ is a rank one local system on the complement
$\tilde X-\pi^{-1}(D)$. However, its monodromy around some components
of $\pi^{-1}(D)$ may be trivial. Denote by $\tilde D\subset\pi^{-1}(D)\subset\tilde X$
the maximal divisor in $\tilde X$ such that $\pi^*\LL$ has nontrivial
monodromy around $\tilde D$. The local systems $\pi^*\LL$ and $\pi^*\LL^\vee$ extend to
local systems on $\tilde U = \tilde X- \tilde D$ denoted by $\tilde \LL$ and $\tilde\LL^\vee$, respectively.
For  $x \in X$, denote $\tilde U_x = \pi^{-1}(x) \cap \tilde U$.

\medskip
\noindent
{\bf Definition 1.}
{\it We say that the local system $\LL$ on $U$  satisfies {\it condition A}
with respect to the resolution of singularities
$\pi : \tilde X\to X$,
if for any edge $F$ of $D$ and
 any point $x$ of $F$ that is not contained in any smaller edge,
 we have  $H^\ell(\tilde U_x; \tilde \LL|_{\tilde U_x})=0$ for $\ell>\text{codim} F-1$. Similarly, we say that
the local system $\LL^\vee$ on $U$   satisfies {\it condition A} with respect to the resolution of singularities
$\pi : \tilde X\to X$,
if for any
edge $F$ of $D$ and any point $x\in F$
that is not contained in any smaller edge, we have
$H^\ell(\tilde U_x; \tilde \LL^\vee|_{\tilde U_x})=0$ for $\ell>\text{codim} F-1$.}

\medskip

It is not hard to check that to verify condition A, it suffices to consider only dense edges $F$.

Notice that if the monodromy of $\LL$ lies in $\{z\in \C^\times\ |\ |z|=1\}$ (in other words, $\LL$ is a unitary local system),
then  $\LL$ satisfies condition A
with respect to $\pi$
if and only if  $\LL^\vee$ satisfies condition A with respect to $\pi$.

\medskip
\noindent
{\bf Theorem 1.}
{\it If both local systems $\LL$  and $\LL^\vee$ on the complement $U$ to the hyperplane-like divisor $D$ in $X$
satisfy condition A with respect to a resolution of singularities
$\pi : \tilde X\to X$, then the intersection cohomology $H^*(X;j_{!*}\LL)$ is naturally isomorphic to}
$H^*(\tilde U; \tilde \LL)$.

\begin{proof}
Let $\tilde j:\tilde U\to\tilde X$ be the open embedding. Since $\tilde \LL$ has non-trivial monodromy
around every component of $\tilde D$ and $\tilde D$ has normal crossings,
we have $\tilde j_*\tilde \LL=\tilde j_!\tilde \LL= \tilde j_{!*}\tilde \LL$
as well as $\tilde j_*\tilde \LL^\vee=\tilde j_!\tilde \LL^\vee= \tilde j_{!*}\tilde \LL^\vee$. We also have
$H^*(X;\pi_*\tilde j_*\tilde \LL) = H^*(\tilde U; \tilde\LL)$ and $H^*(X; \pi_*\tilde j_*\tilde \LL^\vee) = H^*(\tilde U; \tilde\LL^\vee)$.
To simplify notation, we write $\pi_*$ for the direct image in the derived category $R\pi_*$. Thus, $\pi_*\tilde j_*\tilde\LL$ is a complex of sheaves
(more precisely, an object of the corresponding derived category), and $H^*(X;\pi_*\tilde j_*\tilde\LL)$ is its hypercohomology.

The restriction of $\pi_*\tilde j_*\tilde \LL$ to $U$ is $\LL$.
The theorem follows from the  next lemma.

\medskip
\noindent
{\bf Lemma.} {\it
Let $\pi:\tilde X\to X$ be a proper holomorphic map of complex manifolds. For $x\in X$, denote  $\tilde i_x:\pi^{-1}(x)\to \tilde X$ the  embedding of the fiber. Let $\tilde\FF$ be a complex of sheaves with constructible cohomology sheaves on $\tilde X$. (In most examples, $\tilde\FF$ is a perverse sheaf, see
\cite{BBD} for the definition.) Set $\tilde\FF^\vee=\Bbb D\tilde\FF$, where $\Bbb D$ is the Verdier duality functor.}

{\it Suppose $\tilde\FF$ satisfies the following condition B: for every  $\ell>0$ there is an analytic subset $X_\ell\subset X$,
$\codim X_\ell =\ell+1$, such that for any $x\in X- X_\ell$ we have
$H^i(\pi^{-1}(x); \tilde i_x^*\tilde \FF)=0$ and
$H^i(\pi^{-1}(x); \tilde i_x^*\tilde \FF^\vee)=0$ for all $i\ge\ell$ or $i<0$.
Then $\pi_*\tilde \FF=\pi_!\tilde \FF$ is a perverse IC-sheaf on $X$.}

\begin{proof} This lemma is a slight generalization of the results of Goresky and MacPherson (\cite[Section~6.2]{GM}) about small maps. The
argument of \cite{GM} applies without change. Indeed, by \cite[Second theorem of Section~6.1]{GM}, we need to verify that the complex $\pi_*\tilde\FF$
is such that $H^i(\pi_*\tilde\FF)=0$ for $i<0$ and $\codim\supp(H^i(\pi_*\tilde\FF))>i$ for $i>0$, and that the same is true for the dual complex
$\Bbb D\pi_*\tilde\FF=\pi_*\tilde\FF^\vee$. But these conditions are clear by base change. (Note that unlike \cite{GM}, we use the non-self-dual normalization:
for instance, a local system on a smooth manifold is an IC-sheaf in our convention, but it requires cohomological shift in the self-dual normalization.)
\end{proof}

In particular, in the settings of the theorem, the lemma applies to $\FF=\tilde j_*\tilde\LL$, because both $\LL$ and $\LL^\vee$ satisfy condition A.
This concludes the proof of the theorem.
\end{proof}

\medskip
\noindent
{\bf Corollary.} {\it If the
 local systems $\LL$  and $\LL^\vee$ satisfy condition A with respect to each of two resolutions of singularities
$\pi : \tilde X\to X$ and $\pi':\tilde X'\to X$, then
$H^*(\tilde U; \tilde \LL)$ and $H^*(\tilde U'; \tilde \LL')$ are canonically isomorphic.}

\medskip
\noindent
{\bf Example 1.}   Let $X$ be the projective space of dimension $k$.
Let $D\subset X$ be the union of hyperplanes and $\LL$ a rank one local system
on $U=X-D$ with nontrivial monodromy around each hyperplane. Assume that $D$ has normal intersections at all edges
except at vertices. Then both $\LL$ and $\LL^\vee$ satisfy condition A with respect to the minimal resolution
of singularities. Indeed, to obtain $\tilde X$ one has to blow-up
dense vertices of $D$. If $x\in D$ is a dense vertex, then $\tilde U_x$ is nonempty only if $\LL$ has trivial monodromy
across $x$. In that case $\tilde U_x$ is an affine variety of
dimension $k-1$ and  $H^\ell(\tilde U_x; \tilde \LL|_{\tilde U_x})=0$,  $H^\ell(\tilde U_x; \tilde \LL^\vee|_{\tilde U_x})=0$
for $\ell\geq k$.

\medskip
\noindent
{\bf Example 2.}  Let $X$ be $\C^3$. Let $\CC$ be the central arrangement of six planes
\bea
&&
H_1: \ x_1-x_3=0,\qquad
H_2: \ x_1+x_3=0,\\
&&
H_3: \ x_2-x_3=0,\qquad
H_4: \ x_2+x_3=0,\\
&&
H_5: \ x_1-x_2=0,\qquad
H_6: \ x_1+x_2=0
\eea
 with weights $
  a_1=a_2,\quad
a_3=a_4, \quad a_5=a_6,\quad
 a_1+a_3+a_5=0$.
 Let $\LL$ be the local system on $U$ with the monodromy $e^{2\pi i a_j}\neq 1$ around $H_j$.
 Let $\pi:\tilde X\to X$ be the minimal resolution of singularities.
For $x=(0,0,0)$, the space $\tilde U_x$
is the projective plane with four blown-up points  and six lines removed.
We have dim $H^3(\tilde U_x,\tilde \LL)=1$.
This weighted arrangement $\CC$ does not satisfies condition A with respect to the minimal resolution of singularities.

\medskip
\noindent
{\bf Remark.}\ If $D\subset \C^k$ is the union of hyperplanes of a central arrangement and
the monodromy of $\LL$ is close to 1, then
the intersection cohomology $H^*(X;j_{!*}\LL)$ was computed in \cite{KV} as the cohomology of
the complex of flag forms of the arrangement.

\bigskip

The following equivariant version of Theorem 1 holds.
Let $G$ be a finite group, $\rho$ an irreducible representation of $G$.
For a representation $M$ denote by $M^\rho\subset M$ the $\rho$-isotypical component.

Let $\pi:\tilde X\to X$ be a resolution of singularities of $D$ as before.
Assume that  $G$
acts on $X$ and
$\tilde X$ so   that the actions preserve $U, \LL, \tilde U, \tilde \LL$ and commute with the map
 $\pi$.  Then $G$ acts on
$H^*(X;j_{!*}\LL)$ and $H^*(\tilde U; \tilde \LL)$.

 For $x\in X$, we denote by $O_x$ the $G$-orbit of $x$.

\medskip
\noindent
{\bf Definition 2.}
{\it We say that the
local system $\LL$  satisfies {\it condition A} with respect to $\rho$ and
 a resolution of singularities $\pi$ if for any edge $F$
 of $D$ and any  point $x\in F$ that is not contained in any smaller edge,
 we have
 $H^\ell(\cup_{y\in O_x} \tilde U_y; \tilde \LL|_{\cup_{y\in O_x}\tilde U_y})^\rho=0$
 for $\ell>\text{codim} F-1$. Similarly, we say that
the local system $\LL^\vee$ on $U$ satisfies {\it condition A} with respect to $\rho$
 and  a resolution of singularities $\pi$ if for any
 edge $F$ of $D$ and any point $x\in F$ that is not contained in any smaller edge,
 we have
 $H^\ell(\cup_{y\in O_x} \tilde U_y; \tilde \LL^\vee|_{\cup_{y\in O_x}\tilde U_y})^\rho=0$
 for $\ell>\codim F-1$.
}
\medskip

Generally speaking, we can no longer consider dense edges only in the equivariant version of condition A.
In principle, in the equivariant case,
for each edge $F$, it suffices to check
 only the generic points $x\in F$, even though the stabilizer might be different for other points.

\medskip
\noindent
{\bf Theorem 2.}
{\it If both local systems $\LL$  and $\LL^\vee$ on  $U$
satisfy condition A with respect $\rho$  and
 a resolution of singularities $\pi$, then the intersection cohomology $H^*(X;j_{!*}\LL)^\rho$ is naturally isomorphic to}
$H^*(\tilde U; \tilde \LL)^\rho$ as $G$-modules.

\begin{proof} Consider the quotient $X/G$, which may be singular.
The quotient map $q:X\to X/G$ is finite; therefore, the derived direct image $q_*(j_{!*}\LL)$ is a
perverse IC-sheaf on $X/G$. The direct image carries a fiber-wise
action of $G$, so we can take the $\rho$-isotypical component, which is a direct summand
$q_*(j_{!*}\LL)^\rho\subset q_*(j_{!*}\LL)$. Thus $q_*(j_{!*}\LL)^\rho$ is itself an IC-sheaf.

Similarly, there is a direct summand $(q_*\pi_*\tilde j_*\tilde \LL)^\rho\subset q_*\pi_*\tilde j_*\tilde \LL$.
 By the argument used in the proof of Theorem 1,
$(q_*\pi_*\tilde j_*\tilde \LL)^\rho$ is a perverse IC-sheaf. It is thus identified with $q_*(j_{!*}\LL)^\rho$.
\end{proof}

\bigskip
The authors thank A. Levin and M. Falk for helping  develop Example 2.

\end{document}